\documentclass[12pt,twoside]{article}
\usepackage{amsmath, amsfonts, amssymb, amsthm,latexsym}
\usepackage{amscd,graphics}
\usepackage[cp1251]{inputenc}
\usepackage[english]{babel}

\newtheorem{theorem}{Theorem}[section]
\newtheorem{lemma}[theorem]{Lemma}
\newtheorem{corollary}[theorem]{Corollary}
\newtheorem{proposition}[theorem]{Proposition}

\newtheorem{remark}{Remark}

\begin{document}

\newcommand{\iiRGH}{{\rm Rep}\,(G,{\mathcal H})}
\newcommand{\iiRGX}{{\rm Rep}_\circ\,(G,X)}
\newcommand{\iiRGdf}{{\rm Rep}\,(G,d,f)}
\newcommand{\iiRGU}{{\rm Rep}\,(G,\coprod)}
\newcommand{\iiRGUC}{{\rm Rep}_\circ\,(G,\coprod)}
\newcommand{\iiRGUB}{{\rm Rep}_\bullet\,(G,\coprod)}
\newcommand{\iiRGdfC}{{\rm Rep}_\circ\,(G,d,f)}
\newcommand{\iiRGdfB}{{\rm Rep}_\bullet\,(G,d,f)}
\newcommand{\iiSC}[1]{\stackrel{\circ}{#1}}
\newcommand{\iiSB}[1]{\stackrel{\bullet}{#1}}
\newcommand{\iiSCB}[1]{\stackrel{\circ\bullet}{#1}}
\newcommand{\iiSBC}[1]{\stackrel{\bullet\circ}{#1}}
\newcommand{\iiSAI}{\frac{1}{\sqrt{\alpha_i}}}
\newcommand{\iiAI}{\frac{1}{\alpha_i}}



\begin{center}
\LARGE \bf
Locally-Scalar\\
 Representations of Graphs\\
 in the Category of Hilbert Spaces
\end{center}

\begin{center}
\Large\bfseries S.A. Kruglyak~$^*$, A.V. Roiter~$^{**}$
\end{center}

{ \footnotesize
$^*$~Institute of Mathematics of National Academy of Sciences of Ukraine,\\
Tereshchenkovska str., 3, Kiev, Ukraine, ind. 01601
\medskip

$^{**}$~Institute of Mathematics of National Academy of Sciences of Ukraine,\\
Tereshchenkovska str., 3, Kiev, Ukraine, ind. 01601\\
E-mail: roiter@imath.kiev.ua}

\bigskip

\setcounter{section}{0}
\section*{Preface}


In article [1] authors wrote: ``Recently it became evident that
a number of problems of linear algebra allows common formulating and in this
formulating common effective methods of investigations of such problems
appear. It is interesting that these methods appear to be connected with
such concepts as Coxeter--Weyl group and Dynkin schemes.''

One of these problems was a problem of representations of quivers [2] (see
also [3--5]). According to the Gabriel theorem a connected quiver has
finite representative type if and only if corresponding
nonoriented graph is a Dynkin scheme. Coxeter functors defined in
[1] served to explain this fact (which was obtained by P.~Gabriel a
posteriori).

Furher the theory of representations of quivers (in the category of
finite-dimensional vector spaces) was widely developed [6--10].

One may consider representations of quivers also in
metric and, in particular, Hilbert spaces ([11--13]).
But if we use suggesting themselves natural definitions, the problem
of classification of representations of qiuvers (excluding a quiver each
connected component of which consists of one vertex or of two vertices and
and arrow between them) in such categories will become ``wild'' (in terms
of [14], i.~e. including a problem of unitary classification of arbitrary
linear operator or a pair of self-adjoint operators).

In the following we suggest a limitation of local scalarity to the
representations of a graph (or a quiver) in the category of Hilbert spaces,
and after that a theory close to [1--10] arise, though specific for Hilbert
spaces aspects arise.

Note that in some particular case a problem of classification of
locally-scalar representations in other terms studied in [16--17]
(see~example~3).

In \S\S~2,~3 for the category of locally-scalar representations of graph in
the category of Hilbert spaces functors of even and
odd Coxeter reflections are defined. In \S~4 with their help
an analogue of Gabriel theorem is proved.

\section{Basic definitions}

Let ${\mathcal H}$ will be the category of Hilbert spaces, which
objects are separable Hilbert spaces (finite-dimensional or
infinite-dimensional) and morphisms are bounded operators. To each
$\varphi \in \mathcal{H}(A,B)$ uniquely corresponds $\varphi^* \in
\mathcal{H}(B,A)$.

According to Gabriel [2], we will call oriented graphs
as quivers, leaving the term ``graph'' for nonoriented graphs. In both
cases we will admit loops, multiple arrows and edges.

Thus, a quiver $Q$ consists of two sets $Q_v$ (vertices) and $Q_a$
(arrows) and two maps $t$ and $h$ from $Q_a$ to $Q_v$, associated each
arrow $\alpha \in Q_a$ with its tail $t(\alpha)$ and head $h(\alpha)$.

Graph $G$ consists of sets $G_v$ (vertices) and $G_e$ (edges) and
a map $\varepsilon$ from $G_e$ to set of one- and two-element
subsets $G_v$. To each quiver $Q$ in a natural way its graph
$G=G(Q)$ corresponds.

Representation $\pi$ of a quiver $Q$ in an arbitrary category $K$
associates to each vertex $a \in Q_v$ an object $\pi(a) \in K$ and to each
arrow $\alpha \in Q_a$ a morphism $\pi(\alpha) \in K(A,B),$ where
$A=\pi(t(\alpha))$, $B=\pi(h(\alpha))$.

In [1--10] $K={\rm mod}\,k$, where $k$ is a field. Representations of
a quiver $Q$ over field $k$ form a category ${\rm Rep}\,(Q,k)$.

If $\pi$ is a representation of $Q$ in ${\mathcal H}$ then to each arrow
$\alpha : a \to b$, besides operator $\pi(\alpha) : A \to B$, an
operator $\pi^*(\alpha) : B \to A$ corresponds. Thus, in the category
${\mathcal H}$ it is naturally to define a representation $\pi$ of a graph
$G=G(Q)$ which associates to each arrow $a \in G_v$ an object $\pi(a)=H_a
\in {\mathcal H}$, to each edge $\gamma \in G_e$
with $\{i,j\}=\varepsilon(\gamma)$ when $i \neq j$ a pair of interadjoint
linear operators $\pi(\gamma) = \{\Gamma_{ij},\Gamma_{ji}\}$ where
$\Gamma_{ij} : H_j \to H_i$, and when $\varepsilon(\gamma)=\{i\}$ a pair of
interadjoint operators $\pi(\gamma)=\{\Gamma_{ii},\Gamma^*_{ii}\}$.

Consequently, we can identify representations of a quiver $Q$ and of a graph
$G=G(Q)$ in the category ${\mathcal H}$. We will prefer in this paper
to consider in Hilbert spaces representations of graphs, but all
results can be naturally restated for quivers.

Representation $\pi$ of a graph $G$ is said to be {\it finite-dimensional},
if $\dim \pi(i) < \infty$ for all $i \in G_v$.

It is known that graph $G_0$, consisting of one loop, has
infinite-dimensional indecomposable representations (in ${\rm
Rep}\,(G_0,{\mathcal H}))$.

We will call representation $\pi$ of a graph $G$ {\it discrete}, if it
decompose to direct sum (finite or infinite) of finite-dimensional
representations.

Let us construct a category $\iiRGH$ of representations $G$ in ${\mathcal
H}$.

Morphism $C: \pi \to \tilde\pi$ is a family $\{C_i\}_{i\in
G_v}$ of operators $C_i:\pi(i) \to \tilde\pi(i)$ such that diagrams

\[
\begin{CD} H_i   @>\Gamma_{ji}>>  H_j\\ @V{C_i}VV       @VV{C_j}V\\ \tilde
H_i       @>\tilde\Gamma_{ji}>>  \tilde  H_j \end{CD}
\]

\noindent are commutative, i.~e.

\begin{equation}{\label{1}}
C_j\Gamma_{ji}=\tilde\Gamma_{ji}C_i.
\end{equation}

It can be shown that operators $C_i$ realizing an equivalence of two
representations can be chosen as unitary (see, for instance, [15]), i.~e.
equivalent objects of the category $\iiRGH$ are unitary equivalent.

{\it Support} $G^\pi_v = G^\pi$ of a representation $\pi$ is a set
$\{i\in G_v\; |\;\pi(i)\neq 0\}$.

Denote $M_i$ a set of vertices connected with vertex $i$ by edge;
$\overline{M}_i = \{\gamma \in G_e\,|\,i \in \varepsilon(\gamma)\}$; if
$X \subset G_v$ then $M(X) = (\bigcup\limits_{i \in X} M_i) \setminus X$.

Representation $\pi$ {\it faithful} if $G^\pi = G_v$.

Let $\pi \in \iiRGH$. Let $A(\gamma,i) = \Gamma^*_{ji}\Gamma_{ji}$ when
$\varepsilon(\gamma) = \{i,j\}$, and $A(\gamma,i) = \Gamma^*_{ii}\Gamma_{ii} +
\Gamma_{ii}\Gamma^*_{ii}$ when $\varepsilon(\gamma) = \{i\}$, and

$$A_i = \sum\limits_{\gamma \in \overline{M}_i}A(\gamma,i).$$

$A_i$ is a selfadjoint operator in space $H_i = \pi(i)$.
If $\overline{M}_i = \emptyset$ we will consider $A_i = 0$. Note
that if $G$ does not consist loops and multiple edges than

$$A_i = \sum\limits_{j \in M_i}\Gamma_{ij}\Gamma_{ji}.$$

A representation $\pi$ will be named as {\it locally-scalar}
if all operators $A_i$ are scalar, $A_i=\alpha_i I_{H_i}$, where
$I_{H_i}$ is identical operator in space $H_i$; since $A_i$
are positive operators, $\alpha_i \geq 0$.

If $\pi(j)=0$ for all $j \in M_i$ and $\pi(i)\neq 0$, then, obviously,
$\alpha_i=0$. We denote ${\rm Rep}\,(G)$ full subcategory in $\iiRGH$,
which objects are locally-scalar.

Denote by $V_G$ a linear real space, which consists of the
collections $x=(x_i)$ of real numbers $x_i$ ($i \in G_v$);
elements $x$ from $G_v$ we will call $G$-vectors. Vector $x = (x_i)$
we will call {\it positive} ($x > 0$) if $x \neq 0$ and $x_i \geq
0$ for $i \in G_v$. Let us denote $V_G^+ = \{x \in V_G\,|\, x > 0\}$.
Any function $f$ on $G_v$ with the real values can be
identified with corresponding $G$-vector $(f(i))_{i \in G_v}$. $Z_G =
\{d : G_v \to \mathbb{N}_0\}$, $Z_G^+ = Z_G \bigcap V_G^+$.

$G$-vector $d(i)=\dim \pi(i)$ is {\it dimension} of a finite-dimensional
representation $\pi \in \iiRGH$; if $A_i = f(i)I_{H_i}$ for $i \in G_v$,
then $f(i)$ we will call as {\it character} of locally-scalar
representation $\pi$, and $\pi$ will be $f$-representation. Character
is determined uniquely on the support $G^\pi$ of the representation $\pi$
(and ambiguously determined outside of the support). If $G_v = G^\pi$ then
character of the representation is determined uniquely and is
denoted by $f_{\pi}$. In a common case let us denote by
$\{f_\pi\}$ the set of characters of the representation $\pi$.

For graphs there may exist both infinite-dimensional and
finite-dimensional indecomposable representations and the dimensions
of the latter may not be bounded. At the same time, even in the elementary
cases, there exist infinitely many indecomposable representations in the
fixed dimension corresponding to different characters.

We will say that $G$ is (locally-scalar) {\it finitely representable}
in $\mathcal H$, if all its locally-scalar representation are discrete,
dimensions of its indecomposable locally-scalar representations
are bounded in the whole and in each dimension the number
of indecomposable representations with given character is finite.

In what follows we will prove that connected finite graph is finitely
representable if and only if it is Dynkin scheme ($A_n$, $D_n$,
$E_6$, $E_7$, $E_8$), and in this case its indecomposable
representations are uniquely determined by dimension and a value
of character defined on the support of the representation.

\smallskip
{\bf Example 1.} An arbitrary graph $G$ has a locally-scalar
representation $\pi$ in the dimension $d(i) \equiv 1$. Let $\pi(i) =
{\mathbb C}e_i$, $(e_i,e_i) = 1$ and $\pi(\gamma_{ij}) =
\{\Gamma_{ji},\Gamma_{ij}\}$, where $\Gamma_{ji}e_i = e_j$ and
$\Gamma_{ij}e_j = e_i$ ($\varepsilon(\gamma_{ij}) = \{i,j\}$). Then
character $f(i) = |\overline{M}_i|$. $\pi$ indecomposable if and only if
$G$ is connected.

\smallskip
{\bf Example 2.} $G_0$ consists of one vertex $a$ and one loop
$\gamma$. Let $\pi$ be a locally-scalar representation of a graph
$G_0$, for which $\pi(\gamma) = \{\Gamma,\Gamma^*\}$, $f(a) = \alpha$,
i.~e. $\Gamma\Gamma^* + \Gamma^*\Gamma = \alpha I_H$. In [14] it
was shown that all indecomposable locally-scalar
representations no more than two-dimensional; given fixed
positive $\alpha$ in the dimension 2 indecomposable representations
depend on two continuous parameters, so graph $G_0$ is not finitely
representable.

Let us instance nondiscrete representation of the graph $G_0$
which is decomposable, though it can not be decomposed
into direct sum of the indecomposable representations (but only into their
integral). Consider a Hilbert space $H$ with (orthonormal) basis
$\{e_i\}$, $i \in {\mathbb Z}$, and an operator $\Gamma(e_i) =
e_{i+1}$ in it ($\Gamma^*\Gamma + \Gamma\Gamma^* = 2I_H$). Hence $H$ does
not contain finite-dimensional invariant subspaces.

\smallskip
{\bf Example 3.} A problem of (unitary) classification of
locally-scalar representations of the graph $G_n$

\smallskip
\centerline{\includegraphics{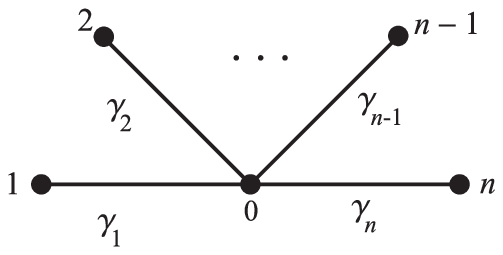}}
\smallskip

\noindent was considered in fact (in other terms) in the works [16--19]. In
[16--17] it was studied a problem of classification up to
unitary equivalence of collections of orthoprojectors $P_1,P_2,\ldots,P_n$
(in the separable Hilbert space $H$) such that
$\sum\limits_{k=1}^n P_k=\alpha I_H$ and, in particular, a problem
of description of the set $\Sigma_n$ of that real numbers
$\alpha$, for which there exists at least one collection of such
orthoprojectors (in the nonzero space).

Corresponding to the orthoprojector $P_i$ a space $H_i={\rm Im}P_i$ and
a natural enclosure (isometry) $\Gamma_i: H_i \to H$ we will get a
locally-scalar representation $\pi$ of the graph $G_n$ if we put
$\pi(i)=H_i$, $i \in \overline{1,n}$ $\pi(0)=H$,
$\pi(\gamma_{i})=\{\Gamma_i,\Gamma_i^*\}$. At that the character $f(i)=1$
at $i\neq 0$ and $f(0)=\alpha$ ($\Gamma_i^*\Gamma_i = I_{H_i}$, $P_i =
\Gamma_i\Gamma_i^*$ and $\sum\limits_{k=1}^n \Gamma_k\Gamma_k^* = \alpha
I_H$).

In [18--19] for representations of the graph $G_n$ were
constructed functors which structure and role during the
description of representations is the same as one of Coxeter functors in
[1]. Using these functors, in [17] full description of the set
$\Sigma_n$ and new results concerning the collections of orthoprojectors
were obtained.

As it was mentioned above, the graph $G_0$ (a loop) has
infinite-dimensional indecomposable representations, but they are not
locally-scalar. As it follows from [17], the graph $G_n$ has
infinite-dimensional indecomposable locally-scalar representations
if and only if $n>4$.

\section{Functors of reflections}

Example~2 implies if $G$ contains a cycle then $G$ is not
finitely representable in $\mathcal{H}$. We will study a question of
finite representability and thus we restict ourselves to graphs not
containing cycles: everywhere below $G$ is a finite connected graph
without cycles (a wood).

Let us fix a decomposition of the set $G_v$ as $\iiSC{G} \coprod
\iiSB{G}$ (univocal up to permutation $\iiSC{G}$ and
$\iiSB{G}$) such that for each $\alpha \in G_e$ one of the vertices from
$\varepsilon(\alpha)$ is situated in $\iiSC{G}$ and the other in
$\iiSB{G}$. Vertexes of the set $\iiSC{G}$ is said to be even and of
$\iiSB{G}$~--- odd.

For each vertex $i\in G_v$ denote as $\sigma_i$ a linear
transformation in the space $V_G$ which is defined by formulas
$(\sigma_ix)_j=x_j$ with $j\neq i$, $(\sigma_ix)_i=-x_i+\sum\limits_{j\in
M_i}x_j$. We will call $\sigma_i$ {\it a reflection} in the vertex $i$.
$W$ is a group of transformations of the space $V_G$ generated by
the reflections $\sigma_i$.

Let us fix a numeration of the vertices of the graph $G$, numbering
first odd and then even vertices. The product of the reflections in all
odd vertices we denote as $\iiSB{c}$, and in the even~--- as $\iiSC{c}$
(remark that reflections in odd (even) vertices commute).
{\it Coxeter transformation on $V_G$} is $c =
\iiSC{c}\iiSB{c}$; $c^{-1} = \iiSB{c}\iiSC{c}$. $\iiSB{c}$ ($\iiSC{c}$)
are said to be odd (even) Coxeter transformations.

Let us provide a denotation for the composition of Coxeter transformations:
$\iiSB{c}_k = \underbrace{\cdots \iiSC{c}\iiSB{c}}_{k \mbox{ times }}$,
$\iiSC{c}_k = \underbrace{\cdots \iiSB{c}\iiSC{c}}_{k \mbox{ times }}$,
$k \in \mathbb{N}$.

Let $X \subset G_v$. ${\rm Rep}_\circ\,(G,X)$~--- full subcategory of
locally-scalar representations $\pi$ in the category ${\rm Rep}\,(G)$ for
which $G^\pi = X$ and the character is positive on $G^\pi \bigcap
\iiSC{G}$. Let $(M(X))^\circ = M(X) \bigcap \iiSC{G}$, $\delta
: (M(X))^\circ \to \mathbb{R}^+$. A case when $(M(X))^\circ = \emptyset$ is
not excluded.

Let us construct a functor

$$\iiSC{F}_{X,\delta} : {\rm Rep}_\circ\,(G,X) \to {\rm Rep}\,(G).$$

In what follows implies that local scalarity of representations is
essential for constructing functor $\iiSC{F}_{X,\delta}$ (if$(M(X))^\circ
= \emptyset$, we will denote it $\iiSC{F}_X$, and if $X = G_v$ then
$\iiSC{F}_G$) and under action of $\iiSC{F}_{X,\delta}$ locally-scalar
representation $\pi$ turn to locally-scalar again.

If Hilbert spaces $H$ and $\widehat H$ are decomposed into orthogonal
sum of subspaces: $H=H^{(1)}\bigoplus H^{(2)}\bigoplus \cdots \bigoplus
H^{(n)}$, $\widehat{H}=\widehat{H}^{(1)}\bigoplus
\widehat{H}^{(2)}\bigoplus \cdots \bigoplus \widehat{H}^{(m)}$ then
arbitrary linear operator $A: H \to \widehat{H}$ can be written
as a matrix $A=[A_{ij}]_{j=\overline{1,n}}^{i=\overline{1,m}}$ where
operators $A_{ij}$ act from the space $H^{(j)}$ to the space
$\widehat{H}^{(i)}$. Such matrices multiplies according to usual rules
of multiplying the block matrices.


Let $\pi \in {\rm Ob}\,\iiRGX$ and $i \in X^\circ \bigcup
(M(X))^\circ$. Let us fix $f \in \{f_\pi\}$ considering $f(j) = \delta(j)$
with $j \in (M(X))^\circ$.

Conditions of local scalarity in the vertex $i$ imply
that operator

$$\Gamma_i^*=[\iiSAI\Gamma_{ig_1},\;
\iiSAI\Gamma_{ig_2},\ldots
\iiSAI\Gamma_{ig_k}],$$

\noindent acting from the space $H^{(i)}=H_{g_1}\bigoplus
H_{g_2}\bigoplus \cdots \bigoplus H_{g_k}$ (the sum is orthogonal,
$\alpha_i = f(i)$, $\{g_1,g_2,\ldots ,g_k\} = M_i$) to the space $H_i$,
has a property $\Gamma_i^*\Gamma_i = I_{H_i}$ (i.~e. operator
$\Gamma_i$ is isometry from the space $H_i$ to the space $H^{(i)}$).
Let $\widehat H_i$~--- orthogonal supplement to ${\rm Im}\,\Gamma_i$.
Then, if $\Delta_i$ is a natural enclosure $\widehat H_i$ to $H^{(i)}$
and

$$\Delta_i^*=[\iiSAI\Delta_{ig_1},\; \iiSAI\Delta_{ig_2},\ldots
\iiSAI\Delta_{ig_k}],$$

\noindent then operator

\begin{equation}{\label{2}}
U_i^*=\begin{bmatrix} {\iiSAI}\Gamma_{ig_1} & {\iiSAI}\Gamma_{ig_2} &
\ldots & {\iiSAI}\Gamma_{ig_k} \\ {\iiSAI}\Delta_{ig_1} &
{\iiSAI}\Delta_{ig_2} & \ldots & {\iiSAI}\Delta_{ig_k} \end{bmatrix}
\end{equation}

\noindent ia a unitary operator, $U_i^*: H^{(i)}\to H_i\bigoplus
\widehat{H}_i$, $U_i^*U_i=I_{H_i\bigoplus \widehat{H}_i}$,
$U_iU_i^*=I_{H^{(i)}}$.

Consider operators $\Delta_{g_t,i}: \widehat{H}_i\to H_{g_t}$ and
$\Delta_{i,g_t}=\Delta^*_{g_t,i}: H_{g_t}\to \widehat{H}_i$.

For $i \in X^\circ \bigcup (M(X))^\circ$ let $\iiSC{\pi}(i) =
\widehat{H}_i$ and $\iiSC{\pi}(j) = H_j$ for
other vertices $j$. $\iiSC{\pi}(\gamma_{ij}) =
\{\Delta_{ji},\Delta_{ij}\}$ for $i \in X^\circ \bigcup (M(X))^\circ$, $j
\in M_i$ and $\iiSC{\pi}(\gamma_{kl}) = \pi(\gamma_{kl})$ for other
edges. As a result we obtain representation $\iiSC{\pi} \in {\rm
Ob}\,\iiRGH$.


Let us show that the representation $\iiSC{\pi}$ will be
locally-scalar.

Since $U_i$ are unitary operators then

\begin{equation}{\label{6}}
  \sum\limits_{j \in M_i} \Delta_{ij}\Delta_{ji}=\alpha_i I_{\widehat{H}_i}
  \mbox{ and } \Gamma_{ji}\Gamma_{ij}+\Delta_{ji}\Delta_{ij}=\alpha_i
I_{H_j} \end{equation}

\noindent with $i \in X^\circ \bigcup (M(X))^\circ,\; j \in M_i$.

Statings (\ref{6}) for vertices $j \in X^\bullet = X \bigcap \iiSB{G}$
imply

\begin{equation}{\label{7}}
  \sum\limits_{i \in M_j} \Delta_{ji}\Delta_{ij} =
  (-\alpha_j+\sum\limits_{i \in M_j} \alpha_i)I_{H_j}
\end{equation}

Thus, if we assume that values of $\iiSC{f}$ matches with values of
$f$ in vertices $i \in X^\circ$, and for $j \in X^\bullet$ $\iiSC{f}(j)$
is equal to the sum of values of $f$ in the vertices adjacent to $j$
minus value of $f$ in the vertex $j$, then $\iiSC{f} \in
\{f_{\iiSC{\pi}}\}$.

Let $\iiSC{F}_{X,\delta}(\pi) = \iiSC{\pi}$. Note that with
different $\delta$ we will obtain different locally-scalar
representation.

Define action of $\iiSC{F}_{X,\delta}$ on the morphisms ${\rm
Rep}_\circ\,(G,X)$. Let $\pi,\;\widetilde{\pi} \in {\rm Ob\,
Rep}_\circ\,(G,X)$. $C = {\{C_k\}}_{k\in G_v}: \pi  \to \widetilde\pi$~---
morphism in ${\rm Rep}_\circ\,(G,X)$. If
$j \not\in X^\circ \bigcup (M(X))^\circ$ let $\iiSC{C}_j = C_j$.


Let $i \in (M(X)^\circ \bigcup X^\circ$ and $\alpha_i = f(i) > 0$,
$\widetilde{\alpha}_i = \widetilde{f}(i) > 0$. Equalities~(\ref{1}) for $j
\in M_i$ imply $\widetilde{\Gamma}_{ij}C_j\Gamma_{ji} =
C_i\Gamma_{ij}\Gamma_{ji}$, therefore
$\sum\limits_{j\in M_i}\widetilde{\Gamma}_{ij}C_j\Gamma{ji} =
(\sum\limits_{j\in M_i}C_i\Gamma_{ij}\Gamma_{ji}) =
\alpha_i C_i$ and, therefore,

\begin{equation}{\label{3}}
C_i=\iiAI\sum\limits_{j \in M_i}\widetilde{\Gamma}_{ij}C_j\Gamma_{ji}.
\end{equation}

On the other hand (\ref{1}) for $j \in M_i$ implies
$\widetilde{\Gamma}_{ij}C_j\Gamma_{ji} =
\widetilde{\Gamma}_{ij}\widetilde{\Gamma}_{ji}C_i$, therefore
$\sum\limits_{j\in M_i}\widetilde{\Gamma}_{ij}C_j\Gamma{ji} =
(\sum\limits_{j\in M_i}\widetilde{\Gamma}_{ij}\widetilde{\Gamma}_{ji})C_i =
\widetilde{\alpha}_i C_i$ and

\begin{equation}{\label{3_}}
C_i=\frac{1}{\widetilde{\alpha}_i}\sum\limits_{j \in
M_i}\widetilde{\Gamma}_{ij}C_j\Gamma_{ji}.
\end{equation}

Thus, if $\alpha_i \neq \widetilde{\alpha}_i\; C_i = 0$.

Let if $\alpha_i \neq \widetilde{\alpha}_i\; \iiSC{C}_i = 0$ and if
$\alpha_i = \widetilde{\alpha}_i$

\begin{equation}{\label{4}}
\iiSC{C}_i=
\iiAI\sum\limits_{j\in M_i}\widetilde{\Delta}_{ij}C_j\Delta{ji}.
\end{equation}

Let us show that the diagram is commutative

\begin{equation}{\label{5}}
\begin{CD}
\widehat{H}_i               @>\Delta_{ji}>>      \widehat{H}_j &=H_j\\
 @V\iiSC{C_i}VV                    @VV\iiSC{C_j}=C_jV\\
\widehat{\widetilde{H}}_i       @>\tilde\Delta_{ji}>>
\widehat{\widetilde{H}}_j &=\widetilde{H}_j
\end{CD}
\end{equation}

\noindent where $j\in M_i$, $i\in \iiSC{G}$, i.~e.
$C_j\Delta_{ji}=\widetilde\Delta_{ji}\iiSC{C_i}$.

From (\ref{4}) follows

$$\widetilde{\Delta}_{ji}\iiSC{C_i}=
\widetilde{\Delta}_{ji}\cdot\iiAI\sum\limits_{k\in M_i}
\widetilde{\Delta}_{ik}C_k\Delta_{ki}$$

From the unitarity of the operator $U^*_i$ (written in form (\ref{2}))
and, analogously, of the operator $\widetilde{U}^*_i$ it is follows that

$$\widetilde{\Gamma}_{ji}\widetilde{\Gamma}_{ik}+
\widetilde{\Delta}_{ji}\widetilde{\Delta}_{ik}=0\mbox{ when } j\neq k\,
\mbox{ and } \widetilde{\Gamma}_{ji}\widetilde{\Gamma}_{ij}+
\widetilde{\Delta}_{ji}\widetilde{\Delta}_{ij}=\alpha_i I_{\widetilde{H}_j},
$$

\noindent and, consequently,

$\widetilde{\Delta}_{ji}\iiSC{C_i}=
\iiAI\sum\limits_{\substack{k\in M_i\\k\neq j}}
[-\widetilde{\Gamma}_{ji}\widetilde{\Gamma}_{ik}C_k\Delta_{ki}]+
[I_{\widetilde{H}_j}-\iiAI\widetilde{\Gamma}_{ji}\widetilde{\Gamma}_{ij}]
\cdot C_j\Delta_{ji}=
C_j\Delta_{ji}-\iiAI\sum\limits_{\substack{k\in M_i\\k\neq j}}
\widetilde{\Gamma}_{ji}\widetilde{\Gamma}_{ik}C_k\Delta_{ki}-
\iiAI\widetilde{\Gamma}_{ji}\widetilde{\Gamma}_{ij}C_j\Delta_{ji}.$

Because of the statings (\ref{1}) we obtain that

$\widetilde{\Delta}_{ji}\iiSC{C_i}=
C_j\Delta_{ji}-\iiAI\sum\limits_{\substack{k\in M_i\\k\neq j}}
C_j\Gamma_{ji}\Gamma_{ik}\Delta_{ki}-
\iiAI C_j\Gamma_{ji}\Gamma_{ij}\Delta_{ji}=
C_j\Delta_{ji}-\iiAI C_j\Gamma_{ji}\cdot
\sum\limits_{k\in M_i}\Gamma_{ik}\Delta_{ki}.$

The last sum equals $0$ because of the orthogonality of the block-rows
of the unitary operator (\ref{2}). Therefore,
$\widetilde{\Delta}_{ji}\iiSC{C_i} = C_j\Delta_{ji}$ and the
diagram (\ref{5}) is commutative.

Commutativity of the dual
diagram can be checked analogously, consequently $\iiSC{C} =
{\{\iiSC{C_k}\}}_{k\in G_v}$ is really a morphism of the category
$\iiRGH$.

Let us show that $\iiSC{F}_{X,\delta}$ retains products
of morphisms and identical morphisms.

Let $C = \{C_k\} : \pi \to \widetilde{\pi}$ and $D = \{D_k\} :
\widetilde{\pi} \to \widetilde{\widetilde{\pi}}$ are morphisms in the
category $\iiRGX$.

Commutativity of the diagrams

$$
\begin{CD}
\widehat{H}_i  @>\Delta_{ji}>> \widehat{H}_j @>\Delta_{ij}>>
\widehat{H}_i \\ @V\iiSC{C}_iVV @V C_j VV @V \iiSC{C}_i VV \\
\widehat{\widetilde{H}}_i @>\widetilde\Delta_{ji}>>
\widehat{\widetilde{H}}_j @>\widetilde\Delta_{ij}>>
\widehat{\widetilde{H}}_i \end{CD} $$

$$
\begin{CD}
\widehat{\widetilde{H}}_i  @>\widetilde{\Delta}_{ji}>>
\widehat{\widetilde{H}}_j @>\widetilde{\Delta}_{ij}>>
\widehat{\widetilde{H}}_i \\ @V\iiSC{D}_iVV @V D_j VV @V \iiSC{D}_i VV \\
\widehat{\widetilde{\widetilde{H}}}_i
@>\widetilde{\widetilde{\Delta}}_{ji}>>
\widehat{\widetilde{\widetilde{H}}}_j
@>\widetilde{\widetilde{\Delta}}_{ij}>>
\widehat{\widetilde{\widetilde{H}}}_i \end{CD} $$

\noindent imply commutativity of the diagram

$$
\begin{CD}
\widehat{H}_i  @>>> \widehat{H}_j @>\Delta_{ij}>>   \widehat{H}_i \\
@V\iiSC{D}_i\iiSC{C}_iVV @V D_jC_j VV @V \iiSC{D}_i\iiSC{C}_i VV \\
\widehat{\widetilde{\widetilde{H}}}_i
@>>>
\widehat{\widetilde{\widetilde{H}}}_j
@>\widetilde{\widetilde{\Delta}}_{ij}>>
\widehat{\widetilde{\widetilde{H}}}_i
\end{CD}
$$

Thus, $\iiSC{D}_i\iiSC{C}_i = \iiAI\sum\limits_{j \in
M_i}\widetilde{\widetilde{\Delta}}_{ij}D_jC_j\Delta_{ji} =
\iiSC{(D_iC_i)}$, i.~e. $\iiSC{F}_{X,\delta}(DC) =
\iiSC{F}_{X,\delta}(D) \cdot \iiSC{F}_{X,\delta}(C)$.

If $C = \{I_{H_k}\}_{k \in G_v}$ then with $j \not\in (M(X))^\circ
\bigcup X^\circ$ $\iiSC{C_j} = C_j = I_{H_j}$, and with $i \in
(M(X))^\circ \bigcup X^\circ$ $\iiSC{C_i} = \iiAI\sum\limits_{j \in
M_i}\Delta_{ij}I_{H_j}\Delta_{ji} = I_{H_i}$, i.~e.
$\iiSC{F}_{X,\delta}(\{I_{H_k}\}) = \{I_{H_k}\}$.

So, we have constructed a functor $\iiSC{F}_{X,\delta}$.

Functor

$$\iiSB{F}_{X,\delta} : {\rm Rep}_\bullet\,(G,X) \to {\rm Rep}\,(G).$$

\noindent is constructed in analogous way.


Let now $\pi$~--- finite-dimensional representation of the graph $G$. It
it easy to count that under the action of the functor
$\iiSC{F}_{X,\delta}$ dimensional function $d(i)$ changes in the following
way: representation $\iiSC{F}_{X,\delta}(\pi)$ has dimension
$\iiSC{c}(d)$, the character of representation $\iiSC{F}_G(\pi)$ is
$\iiSB{c}(f)$.

\[
\iiSC{c}(d)(i) = \left\{
\begin{array}{ll}
-d(i)+\sum\limits_{j\in M_i} d(j) & \mbox{ when }i \in \iiSC{G} \\
d(i) & \mbox{ when }i \in \iiSB{G}
\end{array} \right. \
\]

\begin{equation}{\label{8}}
\iiSB{c}(f)(i) = \left\{
\begin{array}{ll}
-f(i)+\sum\limits_{j\in M_i} f(j) & \mbox{ when }i \in \iiSB{G} \\
f(i) & \mbox{ when }i \in \iiSC{G}
\end{array} \right. \
\end{equation}

Analogously the dimension $\iiSB{c}(d)$ of representation
$\iiSB{F}_{X,\delta}(\pi)$ and the character of representation
$\iiSB{F}_G(\pi)$ can be found.

Given $d \in Z^+_G$, $f \in V_G^+$ consider a full subcategory $\iiRGdf$
in ${\rm Rep}\,(G)$, ${\rm Ob}\,\iiRGdf = \{\pi\,|\, d=d(\pi),f \in
\{f_\pi\}\}$. All representations $\pi$ have common support $X_d = G^\pi =
\{i \in G_v\,|\, d(i) \neq 0\}$.

We will consider these categories under condition
$(d,f) \in S = \{(d,f) \in Z_G^+ \times V_G^+,\,|\,
d(i) + f(i) > 0,\; i \in G_v\}$.

Let us mention that the dimension of a representation is determined
uniquely, but the character, in general, is not, so the objects of $\iiRGdf$
can be considered as pairs $(\pi,f)$, $f \in \{f_\pi\}$.

$\widetilde{\rm Rep}\,(G,d,f) \subset \iiRGdf$~--- full subcategory
which objects are pairs $(\pi,f)$ when representations $\pi$ is
indecomposable in $\iiRGH$. ${\rm Rep}_\circ\,(G,d,f) \subset \iiRGdf$
(${\rm Rep}_\bullet\,(G,d,f) \subset \iiRGdf$)~--- full subcategory,
which object are pairs $(\pi,f)$ with $f(i) > 0$ when $i \in X^\circ =
X \bigcap \iiSC{G}$ ($f(i) > 0$ with $i \in X^\bullet = X \bigcap
\iiSB{G}$). $\widetilde{{\rm Rep}}_\circ\,(G,d,f)$ and $\widetilde{\rm
Rep}_\bullet\,(G,d,f)$ are defined in a natural way.

Let $S_\circ = \{(d,f) \in S\, |\, f(i) > 0 \mbox{ when } i \in
\iiSC{X}_d\}$, $S_\bullet = \{(d,f) \in S\, |\, f(i) > 0 \mbox{
when } i \in \iiSB{X}_d\}$

With $(d,f) \in S_\circ$ ($(d,f) \in S_\bullet$) let us construct
a functor of even and odd reflections

$$\iiSC{F_{df}} : {\rm Rep}_\circ\,(G,d,f) \to {\rm
Rep}_\circ\,(G,\iiSC{c}(d),\iiSC{f}_d),$$

\[
\iiSC{f}_d(i) = \left\{
\begin{array}{ll}
\iiSB{c}(f)(i) & \mbox{ when }i \in \iiSB{X}_d \\
f(i) & \mbox{ when }i \not\in \iiSB{X}_d
\end{array} \right. \
\]

$$\iiSB{F_{df}} : {\rm Rep}_\bullet\,(G,d,f) \to {\rm
Rep}_\bullet\,(G,\iiSB{c}(d),\iiSB{f}_d),$$

\[
\iiSB{f}_d(i) = \left\{
\begin{array}{ll}
\iiSC{c}(f)(i) & \mbox{ when }i \in \iiSC{X}_d \\
f(i) & \mbox{ when }i \not\in \iiSC{X}_d
\end{array} \right. \
\]

\noindent let $\iiSC{F}_{df}(\pi,f) = (\iiSC{\pi},\iiSC{f}_d)$ where
$\iiSC{\pi} = \iiSC{F}_{X_d,\delta}(\pi)$ and $\iiSB{F}_{df}(\pi,f) =
\iiSB{F}_{X_d,\delta}(\pi)$.

Note that if $\iiSB{X}_d = \iiSB{G}$ ($\iiSC{X}_d = \iiSC{G}$) then
$\iiSC{f} = \iiSB{c}(f)$ ($\iiSB{f}_d = \iiSC{c}(d)$).

If $\{C_k\}_{k \in G_v}$ is a morphism from the representation $\pi$ to
the representation $\widetilde{\pi}$ then $\iiSC{F}_{df}(\{C_k\}) =
\iiSC{F}_{X_d,\delta}(\{C_k\})$, $\iiSB{F}_{df}(\{C_k\}) =
\iiSB{F}_{X_d,\delta}(\{C_k\})$.

Formulas (\ref{3}), (\ref{3_}) and (\ref{4}) imply that
functor $\iiSC{F_{df}}$ is faithful and full. Besides, it is easy to
see that $\iiSC{F}_{\iiSC{c}(d),\iiSC{f}_d}\iiSC{F}_{df} \cong {\rm
Id}$ and $\iiSC{F}_{df}\iiSC{F}_{\iiSC{c}(d),\iiSC{f}_d} \cong {\rm Id}$
and functor $\iiSC{F}_{df}$ realizes the equivalence of categories.
Analogous statement about $\iiSB{F}_{df}$ is valid too.

Let us define a category $\iiRGU$. Let ${\rm Ob}\,\iiRGU$ =
$\coprod\limits_{(d,f) \in S}{\rm Ob}\,\iiRGdf$, morphisms between objects
from ${\rm Ob\, Rep}\,(G,d_1,f_1)$ and ${\rm Ob\, Rep}\,(G,d_2,f_2)$
match with morphisms in ${\rm Rep}\,(G,d_1,f_1)$ when $(d_1,f_1) =
(d_2,f_2)$ are missing when $(d_1,f_1) \neq (d_2,f_2)$.

Define following full subcategories in $\iiRGU$: $\iiRGUC$ has as a
set of objects $\coprod\limits_{(d,f) \in S_\circ} {\rm Ob}\,\iiRGdf$,
and $\iiRGUB$~--- the set $\coprod\limits_{(d,f) \in S_\bullet} {\rm
Ob}\,\iiRGdf$. Functors $\iiSC{F}_{df}$ ($\iiSB{F}_{df}$) in a natural way
generate a functor $\iiSC{F}$ ($\iiSB{F}$) on the category $\iiRGUC$
($\iiRGUB$) which realizes the equivalence of the category with itself.

So, we proved:

\begin{theorem}{\label{t21}}

Functors of even and odd reflections

$$\iiSC{F} : \iiRGUC \to \iiRGUC$$

$$\iiSB{F} : \iiRGUB \to \iiRGUB,$$

\noindent are defined; they realize the equivalence of the category with
itself and the equivalence of following full subcategories:

$$\iiSC{F} : {\rm Rep}_\circ\,(G,d,f) \to {\rm
Rep}_\circ\,(G,\iiSC{c}(d),\iiSC{f}_d) \mbox{ when } (d,f) \in S_\circ$$

$$\iiSB{F} : {\rm Rep}_\bullet\,(G,d,f) \to {\rm
Rep}_\bullet\,(G,\iiSB{c}(d),\iiSB{f}_d) \mbox{ when } (d,f) \in
S_\bullet.$$

At that $(\iiSC{F})^2 \cong {\rm Id}$, $(\iiSB{F})^2 \cong {\rm Id}$.

\end{theorem}

Let $g \in G_v$, $\Pi_g$~--- simplest representation of a graph $G$:
$\Pi_g(g) = \mathbb{C}$, $\Pi_g(i) = 0$ when $i \neq g$, $i \in
G_v$. The characters of the representations $\Pi_g$ we will denote as
$f_g$: $f_g(g) = 0$, at the same time assuming that $f_g(i) > 0$ when
$i \neq g$.

For the simplest representation $\Pi_g$ the dimension $d_g(g) = 1$,
$d_g(i) = 0$ when $i \neq g$.

Objects ($\Pi_g,f_g$) are said to be simplest objects of the
category $\iiRGU$. If $g \in \iiSB{G}$ ($g \in \iiSC{G}$) then
$(\Pi_g,f_g) \in {\rm Ob}\,\iiRGUC$ ($(\Pi_g,f_g) \in {\rm
Ob}\,\iiRGUB$).

\section{Coxeter tranformations}

Denote as $B$ a quadratic form on the space $V_G$ defined by formula
$B(x)=\sum\limits_{i\in G_v}x_i^2-\sum\limits_{\gamma_{ij}\in G_e}
x_ix_j$, and as $<\;,\;>$~--- corresponding symmetric bilinear form.
Form $B(x)$ is named {\it the Tits form} of the graph $G$.

In propositions~3.1--3.3 we collect well-known results from
[1] which we will use further.

\begin{proposition}{\label{Gel1}} $ $

1. Let $i \in G_v$, then $\sigma_i(x) = x-2<\overline{i},x>\overline{i}$,
$\sigma_i^2 = 1$ .

2. Group $W$, generated by reflections $\sigma_i$, retains the
integral lattice in $V_G$ and retains the quadratic form $B$.

3. If the form $B$ is positive defined then group $W$ is finite.

4. Form $B$ is positive defined for graphs $A_n$, $D_n$, $E_6$, $E_7$,
$E_8$ (Dynkin graphs) and only for them.

\end{proposition}

For each $i \in G_v$ we denote as $\bar{i}$ a vector in $V_G$ such that
$(\bar{i})_j = 0$ when $i \neq j$ and $(\bar{i})_i = 1$.

Vector $x\in V_G$ is said to be a {\it root} if for some
$i\in G_v$ and $w\in W$ we have $x=w \overline{i}$. Vectors $\overline{i}$
are {\it simple roots}. Root $x$ is {\it positive} if $x>0$.

\begin{proposition}{\label{Gel2}} $ $

1. If $x$~--- root then $x$~--- integral vector and $B(x)=1$.

2. If $x$~--- root then $-x$~--- root.

3. If $x$~--- root then either $x>0$ or $(-x)>0$.

\end{proposition}

\begin{proposition}{\label{Gel3}}

If the form $B$ for a graph $G$ is positive defined then:

1) transformation $c$ in the space $V_G$ has no nonzero invariant
vectors;

2) if $x\in V_G,\; x\neq 0$ then for some natural $k$ vector
$c^k x$ is not positive.

\end{proposition}

Let us return to locally-scalar representations of a graph $G$ in the
category of Hilbert spaces.

\begin{lemma}

If $G$ is Dynkin graph and $d = (d_i)$~--- its positive not
simple root then $d$ results from the simple root by the sequence
of even and odd Coxeter transformations.

\end{lemma}

{Proof.} Let $m$~--- the minimal natural number with the property:
$c^m(d) \not\in G_V^+$ (such $m$ can be found by the
proposition~\ref{Gel3}). Applying to $d$ sequentially transformations
of odd and even reflections we will obtain positive root
$\widetilde{d}$ such that the next vector will be negative. Then, as is
well-known, $\widetilde{d}$ will be simple root (transfer from
the positive root to the negative one is made only through the
simple root). Let $\widetilde{d} = \iiSB{c}_k(d)$ (or $\iiSC{c}_k(d)$),
then $d = \iiSB{c}_k^{-1}(\widetilde{d})$ (or
$\iiSC{c}_k^{-1}(\widetilde{d})$) and the lemma is proved.

\medskip


\begin{lemma}\label{l_dir_sum}

Let $\pi \in {\rm Rep}\,(G)$, $i \in G_v$, $d(i) \neq 0$, $f(i) = 0$.
Then $\pi = \pi_i \bigoplus \pi'$, where $G^{\pi_i} = i$ and $G^{\pi'}
\not\ni i$ (a case $\pi' = 0$ is not excluded).

\end{lemma}

Indeed, $f(i) = 0$ implies $\Gamma_{ij} = 0$ when $\gamma_{ij}
\in \overline{M}(i)$, and $d(i) \neq 0$ implies $\pi = \pi_i \bigoplus
\pi'$.

\begin{proposition}{\label{p_dis}}

All locally-scalar representation of a Dynkin graph $G$ in the category
of Hilbert spaces are discrete.

\end{proposition}

{\bf Proof} will be carried out by the induction on $n=|G_v|$.
When $n=1$ the statement is trivial (separable Hilbert
space is orthogonal direct sum, finite or infinite, of the
one-dimensional spaces).

Let the statement is proved for Dynkin graphs with the number of vertices,
which is $\leq n-1$.

Let $n>1$ and $\pi \in {\rm Ob\, Rep}\,(G)$. If $G^\pi \neq G_v$ then it is
possible to turn from $\pi$ to the representation $\widetilde{\pi}$
of the graph $\widetilde{G}$ with the less number of vertices
($\widetilde{G}_v = G^\pi$) and to make use of the presumption of induction.

Let $G^\pi = G_v$. Presume that $f \in \{f_\pi\}$,
$f(g) = 0$, $g \in G_v$, then by the lemma~\ref{l_dir_sum} $\pi = \pi_g
\bigoplus \pi'$ and we can use the presumption of induction again.

Let $G^\pi = G_v$ and the character $f(i) > 0$ when $i \in G_v$. Applying to
$(\pi,f)$ by turns functors $\iiSC{F}_G$, $\iiSB{F}_G$ we
will obtain, subject to (\ref{8}):

a) either $c^k(f)$ is positive for all $k$, which contradicts to
proposition~\ref{Gel3};

b) or after applying functor one of the conditions
$G^{\widetilde{\pi}} = G_v$, $\widetilde{f}(i) > 0$ when $i \in G_v$ will
fail (and in this case we will return to the situation considered earlier).

\medskip


\begin{proposition}\label{p36}

Any locally-scalar indecomposable representation $\pi$ of a Dynkin graph
$G$ results from the some simplest representation $\Pi_{g}$ by the
functors of even and odd reflections. More precisely: if
$(\pi,f) \in \widetilde{\rm Rep}\,(G,\coprod)$ then there exists
a sequence $(\pi^{(k)},f^{(k)})$ of objects of the category $\widetilde{{\rm
Rep}}\,(G,\coprod)$, $k \in \overline{0,n}$ such that $(\pi^{(0)},f^{(0)})
= (\Pi_{g},f_{g})$ (let, for definiteness, $g \in \iiSC{G}$) and
with odd $k$ $\iiSB{F}(\pi^{(k-1)},f^{(k-1)}) = (\pi^{(k)},f^{(k)})$, and
with even $k$ $\iiSC{F}(\pi^{(k-1)},f^{(k-1)}) =
(\pi^{(k)},f^{(k)})$, $k \in \overline{1,n}$, and $(\pi^{(n)},f^{(n)}) =
(\pi,f)$.

\end{proposition}

{\bf Proof.}

Let $\pi$ has a dimension $d$ and a character $f$. Let $n$~--- minimal
number such that $\{\iiSC{c}_n(d),\iiSB{c}_n(d)\} \not\subset V_G^+$.
We will name $n$ as {\it growth} of the locally-scalar
representation $\pi$.

Proposition will be proved by induction on
$n$. Let $n = 1$, then $f(g) = 0$ for some $g \in G^\pi \bigcap
\iiSC{G}$ (in the contrary case we will apply to $\pi$ the functor
$\iiSC{F}_{X,\delta}$, where $X = G^\pi$, and obtain
$d(\iiSC{F}_{X,\delta}(\pi)) = \iiSC{c}(d(\pi)) \in V_G^+$). Then from
indecomposability $\pi$ and the lemma~3.5 it is follows that $G^\pi =
\{g\}$, i.~e. $\pi = \Pi_g$, $(\pi,f)$ is a simplest object.

Let the statement is proved for locally-scalar
representations of the growth $\leq n-1$ and $\pi$ has a growth $n \geq 1$.

By the lemma~\ref{l_dir_sum} $(\pi,f) \in \iiRGUC$ and $(\pi,f) \in
\iiRGUB$ and to the pair $(\pi,f)$ both functors $\iiSC{F}$
and $\iiSB{F}$ may be applied. We will apply that functor, which will entail
a new pair with representation of lower growth. Then we make use of
the presumption of induction.

\medskip

Theorem~\ref{t21} and the proposition~\ref{p36} imply

\begin{corollary}\label{c_equiv}

In the category $\widetilde{\rm Rep}\,(G,d,f)$ all objects are equivalent.

\end{corollary}


Let $\bar{i}$ is a simple root from $V_G$ and, for definiteness, $i \in
\iiSC{G}$. Then for some natural number $l(i)$ the vector
$c^{l(i)}(\bar{i}) > 0$ and $c^{l(i)+1}(\bar{i}) \ngtr 0$.
We will assume that $l(i)$ is a minimal number with such property. As we
mentioned above, the transfer from the positive root $x$ to
nonpositive root $y$ by alternate transformations
$\iiSB{c}$ and $\iiSC{c}$ is made by passing through the simple root.
Minimal nonzero number of steps, for which the simple root $\bar{i}$
transforms to the simple root $\bar{j}$ by transformations $\iiSB{c}$,
$\iiSC{c}$, we define as $h(i)$: $\iiSB{c}_{h(i)}(\bar{i}) = \bar{j}$.
There are two possibilities:

a) $h(i)=2l(i)$. In this case vertex $j$ is odd ($\bullet$). In the series
of roots $\bar{i},\, \iiSB{c}(\bar{i}),\,
\iiSB{c}_2(\bar{i}),\ldots ,\, \iiSB{c}_{2l(i)}(\bar{i}) = \bar{j}$ any
root results from the simple root ($\bar{i}$ or $\bar{j}$) by
transformations $\iiSB{c}$, $\iiSC{c}$ for the number of steps $\leq
\frac{h(i)}{2}$. Middle root in series we will obtain, for definiteness,
from odd simple root (in $\frac{h(i)}{2} = l(i)$ number of steps).

b) $h(i)=2l(i)+1$, $\iiSB{c}_{2l(i)+1}(\bar{i}) = \bar{j}$. In this case
$\bar{j}$ has the same parity as $\bar{i}$ (in particular, it is possible
that $\bar{i}=\bar{j}$). Any root in series $\bar{i},\,
\iiSB{c}(\bar{i}),\ldots ,\, \iiSB{c}_{2l(i)+1}(\bar{i}) = \bar{j}$
results from the simple root ($\bar{i}$ or $\bar{j}$) by transformations
$\iiSB{c}$, $\iiSC{c}$ in $< \frac{h(i)}{2}$ number of steps.


Let $S_G$~--- set of simple roots of a graph $G$, $u_{g}(i) \in
(\mathbb{R}^+)^{G_v \setminus \{g\}}$ (we will construct a character
$f_{g}(i)$ by function $u_{g}$ assuming $f_{g}(i) = u_{g}(i)$ when $i \in
G_v \setminus \{g\}$ and $f_{g}(g)=0$). Let $N_G=\{ (\bar{i},k,u_i)\, |\,
\bar{i}\in S_G;\; k\leq \frac{h(i)}{2} \mbox{ for } i\in \iiSB{G} \mbox{ and
} k<\frac{h(i)}{2}\mbox{ for } i\in \iiSC{G};\; u_i\in (\mathbb{R}^+)^{G_v
\setminus \{i\}} \}$, and ${\rm Ind}\,G$~--- set of indecomposable
locally-scalar representations of a graph $G$ defined up to
unitary equivalence.

By the simple root $\bar{g}$ and function $u_{g}(i)$ by the sequence
of Coxeter reflections in $k$ steps we will obtain from the simplest
object $(\Pi_{g},f_{g})$ an object $(\pi,f)$ of the category
$\widetilde{\rm Rep}\,(G,\coprod)$ ($\pi$ is indecomposable in ${\rm
Rep}\,(G)$ representation from ${\rm Ind}\,G$).

So, we have defined a map

$$\varphi : N_G\to {\rm Ind}\,G.$$

All indecomposable locally-scalar representations of Dynkin graph
are obtained in this way (see proposition \ref{p36}).


Let $\overline{{\rm Ind}}\,G$ is a subset in the set ${\rm Ind}\,G$ of
faithful representations, $T_G$~--- subset of the simple roots $\bar{i}$
and $L_i$~--- set of that values $k$, for which triples $(\bar{i},k,u_i)\in
\varphi^{-1}(\overline{{\rm Ind}}\,G)$.

Following statement holds:

\begin{theorem}{\label{t37}}

The map $\varphi : N_G\to {\rm Ind}\,G$ is surjection, at that each
faithful representation $\pi$ from ${\rm Ind}\,G$ has unique
$\varphi$-inverse image; for nonfaithful $\pi$ vector $\bar{i}$ and the
number $k$ in $\varphi^{-1}(\pi)$ are determined uniquely, and function
$u_i$~--- ambiguously.

\end{theorem}

{\bf Proof} we will obtain by direct count for Dynkin graphs
$E_6$, $E_7$, $E_8$ and by induction for graphs $A_n$ and $D_n$.

\medskip
1. For graph $A_n$

\smallskip
\centerline{\includegraphics{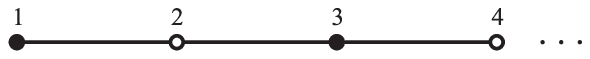}}
\smallskip

$h(i)=n$ for $i=\overline{1,n}$;

when $n=2m+1$, or $n=2m$ with even $m$: $T_G=\{\overline{m+1}\}$,
$L_{m+1}=\{m\}$;

when $n=2m$ with odd $m$: $T_G=\{\overline{m}\}$, $L_m=\{m\}$.

\medskip
2. For graph $D_n\; (n\geq 4)$

\smallskip
\centerline{\includegraphics{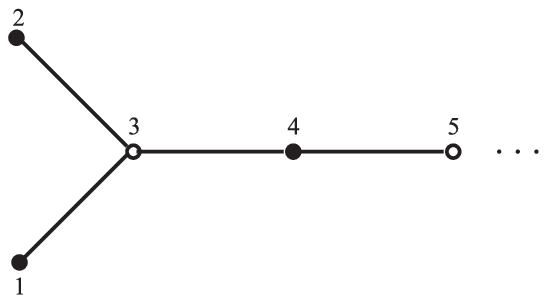}}
\smallskip

$h(i)=2n-3$, $i=\overline{1,n}$;

when $n=2m+1$: $T_G=\{\overline{3},\overline{4},\ldots,\overline{m+2}\}$,
$L_3=\{2m-2,2m-1\}$, $L_4=\{2m-3,2m-2\}$, $\ldots$, $L_{m+1}=\{m,m+1\}$,
$L_{m+2}=\{m\}$;

when $n=2m$: $T_G=\{\overline{3},\overline{4},\ldots,\overline{m+1}\}$,
$L_3=\{2m-3,2m-2\}$, $\ldots$, $L_{m+1}=\{m-1,m\}$.

\medskip
3. For graph $E_6$

\smallskip
\centerline{\includegraphics{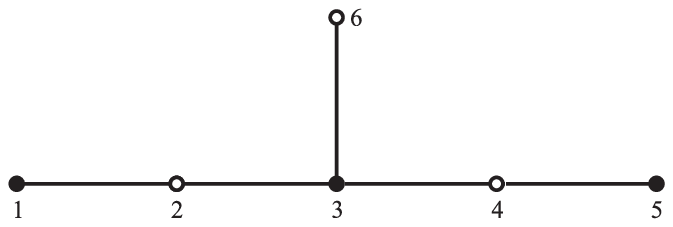}}
\smallskip

$h(i)=11$, $i=\overline{1,6}$;
$T_G=\{\overline{2},\overline{3},\overline{4},\overline{6}\}$,
$L_2=\{5\}$, $L_3=\{2,3,4,5\}$, $L_4=\{5\}$, $L_6=\{4\}$.

\medskip
4. For graph $E_7$

\smallskip
\centerline{\includegraphics{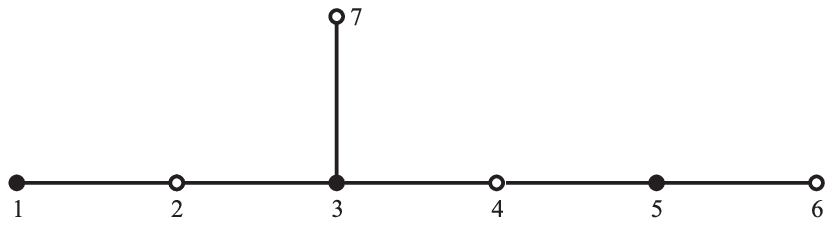}}
\smallskip

$h(i)=17$, $i=\overline{1,7}$;
$T_G=\{\overline{2},\overline{3},\overline{4},\overline{5}, \overline{7}\}$,
$L_2=\{5,6,7\}$, $L_3=\{3,4,5,6,7,8\}$, $L_4=\{3,6,7,8\}$, $L_5=\{7\}$,
$L_7=\{4,8\}$.

\medskip
5. For graph $E_8$

\smallskip
\centerline{\includegraphics{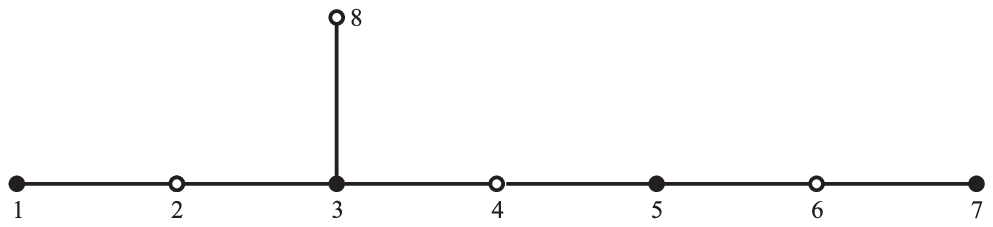}}
\smallskip

$h(i)=29$, $i=\overline{1,8}$;
$T_G=\{\overline{1}, \overline{2},\overline{3},\overline{4},\overline{5},
\overline{6}, \overline{8}\}$, $L_1=\{7,12,13\}$,
$L_2=\{5-8,11-14\}$, $L_3=\{4-14\}$, $L_4=\{3,4,7-14\}$,
$L_5=\{8-13\}$, $L_6=\{11,12\}$, $L_8=\{9,10,13,14\}$.

\medskip

Hence, at the time of proving the theorem {\ref{t37}}, we
proved

\begin{corollary}

For any Dynkin graph $G$ there exists one-to-one
correspondence between faithful indecomposable locally-scalar
representations of graph $G$ and triples $(\bar{i},k,u_i)\in N_G$
where $\bar{i}\in T_G,\; k\in L_i$.

\end{corollary}


Consider in space $V_G$ a scalar product $xy$, meaning
vectors $\bar{i}\; {(i\in G_v)}$ to form an orthonormal basis in $V_G$.
So, if $x = (x_i),\; y = (y_i)$ then $xy = \sum\limits_{i\in
G_v}x_iy_i$, $x = \sum\limits_{i\in G_v}x_i\bar{i} = \sum\limits_{a\in
\iiSB{G}}x_a\bar{a}+\sum\limits_{b\in \iiSC{G}}x_b\bar{b}$.

Let us provide denotations $x^\bullet = \sum\limits_{a\in
\iiSB{G}}x_a\bar{a}$, $x^\circ = \sum\limits_{b\in \iiSC{G}}x_b\bar{b}$ \
($x=x^\bullet + x^\circ$).

Making allowed item 1) of the proposition \ref{Gel1} it is easy to
obtain that $\iiSB{c}(x) = x-2\sum\limits_{a\in
\iiSB{G}}{<\bar{a},x>}\bar{a} = x-2\sum\limits_{a\in
\iiSB{G}}{<\bar{a},x^\bullet + x^\circ >}\bar{a} = x-2\sum\limits_{a\in
\iiSB{G}}{<\bar{a},x^\bullet >}\bar{a} - 2\sum\limits_{a\in
\iiSB{G}}{<\bar{a},x^\circ >}\bar{a} = x-2x^\bullet - 2\sum\limits_{a\in
\iiSB{G}}{<\bar{a},x^\circ >}\bar{a} = x^\circ -x^\bullet
-2\sum\limits_{a\in \iiSB{G}}{<\bar{a},x^\circ >}\bar{a}$. (Here we
uses that fact that $<\bar{a},\bar{a}>=1$ and $<\bar{a},\bar{a}'>=0$
when $a\neq a'$, $a,\; a' \in \iiSB{G}$).

Finally,
$$\iiSB{c}(x) = x^\circ -x^\bullet - 2\sum\limits_{a\in
\iiSB{G}}<\bar{a},x^\circ >\bar{a}.$$

Analogously,
\begin{equation}{\label{12}}
\iiSC{c}(x)=x^\bullet - x^\circ - 2\sum\limits_{b\in
\iiSC{G}}<\bar{b},x^\bullet
> \bar{b}. \end{equation}

Let us prove

\begin{proposition}

1. If vectors $x,y\in V_G$ and $x^\bullet y^\bullet =x^\circ y^\circ$ then
${[\iiSB{c}(x)]^\bullet}\cdot {[\iiSC{c}(y)]^\bullet} =
{[\iiSB{c}(x)]^\circ}\cdot {[\iiSC{c}(y)]^\circ}$;

2. If $\pi$~--- finite-dimensional $f$-representation of a graph $G$
with dimension $d(i),\; i\in G_v$ then $f^\bullet \cdot d^\bullet =
f^\circ \cdot d^\circ$

\end{proposition}

{\bf Proof.}

1. Equalities~({\ref{12}}) imply

$$[\iiSB{c}(x)]^\bullet\cdot [\iiSC{c}(y)]^\bullet =
(-x^\bullet -2\sum\limits_{a\in \iiSB{G}}{<\bar{a},x^\circ >}\bar{a})
y^\bullet = -x^\bullet y^\bullet -2\sum\limits_{a\in
\iiSB{G}}{<\bar{a},x^\circ >}y_a^\bullet = -x^\bullet y^\bullet
-2{<y^\bullet ,x^\circ >},$$

$$[\iiSB{c}(x)]^\circ \cdot
[\iiSC{c}(y)]^\circ = x^\circ \cdot (-y^\circ -2\sum\limits_{b\in
\iiSC{G}}{<\bar{b},y^\bullet
>}\bar{b}) = -x^\circ y^\circ -2\sum\limits_{b\in
\iiSC{G}}{<\bar{b},y^\bullet
>}x_b^\circ = -x^\circ y^\circ -2{<x^\circ ,y^\bullet >},$$

\noindent and since
$x^\bullet y^\bullet = x^\circ y^\circ $ and $<y^\bullet ,x^\circ >=<x^\circ
,y^\bullet >$ then $[\iiSB{c}(x)]^\bullet\cdot [\iiSC{c}(y)]^\bullet =
[\iiSB{c}(x)]^\circ \cdot [\iiSC{c}(y)]^\circ $.

2. The condition of the local scalarity $\sum\limits_{j\in
M_i}\Gamma_{ij}\Gamma_{ji} = \alpha_i I_{H_i}$ for finite-dimensional
representation implies $\sum\limits_{j\in
M_i} {\rm Tr}(\Gamma_{ij}\Gamma_{ji}) = \alpha_id_i$ ($\alpha_i=f(i),\;
d_i=d(i),\; i\in G_v$).

Summig such equalities at first by odd vertices and then
by even ones, we will obtain

$$f^\bullet d^\bullet =
\sum\limits_{\{\gamma_{ij}\gamma_{ji}\}\subset G_e} {\rm Tr}
(\Gamma_{ij}\Gamma_{ji}) =
\sum\limits_{i \in G_v} \alpha_i d_i
= f^\circ d^\circ .$$

\section{Analogue of the Gabriel theorem}

\begin{proposition}{\label{p41}}

If a graph $G$ is one of the extended Dynkin graphs $\widetilde{A}_n$,
$\widetilde{D}_n$, $\widetilde{E}_6$, $\widetilde{E}_7$, $\widetilde{E}_8$
then dimensions of the indecomposable locally-scalar representations  of
the graph $G$ are not bounded in a whole, so that $G$ is not finitely
representable in $\mathcal{H}$.

\end{proposition}

{\bf Proof.}

From the root theory of extended Dynkin graphs it is follows that
if $g \in \iiSC{G}$ then when $t > 1$ $c^t(\bar{g}) \in V_G^+$ and it is
not a simple root (and if $g \in \iiSB{G}$ then when $t > 1$
$c^{-t}(\bar{g}) \in V_G^+$ and it is not a simple root) (see also
[21]), hence $c^{t_1}(\bar{g}) \neq c^{t_2}(\bar{g})$ when $t_1 \neq t_2$.

Then if we fix $g \in \iiSC{G}$ and simplest object $(\Pi,f) \in
\iiRGU$ we will obtain an infinite sequence of objects
$\cdots \iiSC{F}\iiSB{F}(\Pi,f)$ in which dimensions of representations
are different.

\medskip


Let us point to a connection between indecomposable representations
of quivers which graph is a Dynkin scheme and locally-scalar
representations of these graphs.

Let $\pi^\Lambda$ is a finite-dimensional representation of a quiver $Q$
over the field $\mathbb{C}$. If spaces of the representation are Hilbert
(i.~e. unitary, since they are finite-dimensional) then the
representation $\pi^\Lambda$ can be naturally continued to the
representation $\pi$ of the graph $G = G(Q)$ in the category of Hilbert
spaces.

If $\pi^\Lambda$ is equivalent
in ${\rm Rep}\,(Q,\mathbb C)$ to the representation
$\widetilde{\pi}^\Lambda$ such that the continued representation
$\widetilde{\pi}$ of a graph $G$ is locally-scalar, we will say that
the representation $\pi^\Lambda$ of a quiver $Q$ is {\it unitarizable}.

\begin{proposition}{\label{p42}}

Graph $G$ is a Dynkin graph if and only if with $G = G(Q)$
any finite-dimensional indecomposable representation $\pi^\Lambda$
of a quiver $Q$ is unitarizable.

\end{proposition}

{\bf Proof}. If a graph $G$ is a Dynkin graph then dimension
$d(i)$ of the representation $\pi^{\Lambda}$ is a root which can be
obtained from a some simple root $\bar{g}$ by the transformations
$\iiSB{c}$ and $\iiSC{c}$. From the simplest representation $\Pi_{g}$ in
this dimension by functors $\iiSB{F}$ and $\iiSC{F}$
a locally-scalar representation $\widetilde{\pi}$ can be obtained.
Representations $\pi^{\Lambda}$ and $\widetilde{\pi}^{\Lambda}$ of a
quiver $Q$ are equivalent in ${\rm Rep}\,(Q,\mathbb C)$ since for Dynkin
graphs in the dimension $d(i)$ there exist a unique up to the
equivalence indecomposable representation ([1,2]).

The condition of unitarizability can be restated also in such a way: for
the representation $\pi^{\Lambda}$ of a quiver $Q$ it is possible to
define a scalar product in the representation spaces so as the continued
representation of a graph $G$ in the category of Hilbert spaces will be
locally-scalar.

We will prove the second part of the statement of
the proposition~\ref{p42} if we show that for any extended
Dynkin graph $G$ (and accordingly for $G' \supset G$) and any quiver $Q$,
such that $G(Q)=G$, there exists indecomposable representation
$\pi^\Lambda$ of the quiver $Q$ which will not be unitarizable.

If $(d(i))_{i \in G_v}$~--- minimal imaginary root of the graph $G$ (i.~e.
minimal positive root of the equation $B_G(x)=0$) then such
representation exists in the dimension $2d(i)$.

Let, for definiteness, orientation of $Q$ will be such that for any
even vertex $i$ all incident with $i$ arrows are ``come into'' vertex,
and for odd~--- ``come out'' it.

If all incident with $i$ arrows ``come into'' vertex $i$ then there exists
(see, for instance, [10]) a representation $\pi^\Lambda$, for which

\begin{equation}{\label{13}}
\pi^\Lambda(\gamma_{ij}) =
\begin{bmatrix}
\pi_1^\Lambda(\gamma_{ij}) & \multicolumn{1}{|c}{X(\gamma_{ij})} \\
\hline
0                          & \multicolumn{1}{|c}{\pi_2^\Lambda(\gamma_{ij})}
\end{bmatrix},
\end{equation}

\noindent $\pi^\Lambda(i) = H_1^{(i)}\bigoplus H_2^{(i)}$,
$\pi^\Lambda(j) = H_1^{(j)}\bigoplus H_2^{(j)}$,
$\pi_1^\Lambda(\gamma_{ij}): H_1^{(j)}\to H_1^{(i)}$,
$\pi_2^\Lambda(\gamma_{ij}): H_2^{(j)}\to H_2^{(i)}$,
$d_i = {\rm dim}\, H_1^{(i)} = {\rm dim}\, H_2^{(i)}$ for $i\in G_v$.

Assume that representation $\pi^\Lambda$ is unitarizable, i.~e. in
the spaces of representation we can define a scalar product so as
the continued representation $\pi$ of the graph $G$ will become a
locally-scalar.

Selecting an orthogonal basis in spaces $H_1^{(i)}$ and supplementing it to
an orthogonal basis in $H_1^{(i)}\bigoplus H_2^{(i)}$ we will obtain new
matrices of the representation $\widetilde{\pi}^\Lambda$, remaining
reduced and satisfying to the condition of the local scalarity. Saving
generality, we will consider that the representation
$\pi^\Lambda$ itself has this property.

Let $a_{kl}^{(ij)}$~--- matric elements of the matrix
$\pi_1^\Lambda(\gamma_{ij})$, hence sum $\sum\limits_{kl}
\bar{a}_{kl}^{(ij)} a_{kl}^{(ij)}$ is a sum of squares of the column lengths
of the matrix $\pi_1^\Lambda(\gamma_{ij})$ and, therefore,
$$A\equiv \sum\limits_{\gamma_{ij}} \sum\limits_{kl}
\bar{a}_{kl}^{(ij)} a_{kl}^{(ij)} = \sum\limits_{j\in
\iiSB{G}}d_j\alpha_j.$$

If $x_{kl}^{(ij)}$ and $b_{kl}^{(ij)}$ are matric elements of the matrices
$X(\gamma_{ij})$ and, respectively, $\pi_2^\Lambda(\gamma_{ij})$ then for
residuary columns of the matrices of the representation we have

$$C+B\equiv
\sum\limits_{\gamma_{ij}} \sum\limits_{kl} \bar{x}_{kl}^{(ij)}
x_{kl}^{(ij)} + \sum\limits_{\gamma_{ij}} \sum\limits_{kl}
\bar{b}_{kl}^{(ij)} b_{kl}^{(ij)} =
\sum\limits_{j\in \iiSB{G}}d_j\alpha_j,$$

\noindent where $C > 0$; thus $A > B$.

Analogously actions with rows (summing by even vertices) gives us
an inequality $B>A$. Thus, asuumption about local scalarity
of the representation $\pi$ leads us to contradiction.


Above we defined the sets $S$, $S_\circ$, $S_\bullet$. Let us define
a map $\iiSC{\mathcal{F}} : S_\circ \to S_\circ$ for which
${\iiSC{\mathcal{F}}(d,f) = (\iiSC{c}(d),\iiSC{f}_d)}$, and a map
$\iiSB{\mathcal{F}} : S_\bullet \to S_\bullet$ for which
$\iiSB{\mathcal{F}}(d,f) = (\iiSB{c}(d),\iiSB{f}_d)$. If the pair
$\underbrace{\cdots\iiSC{\mathcal{F}}\iiSB{\mathcal{F}}}_{k \mbox{ times
}}(d,f)$
is defined then we will denote it as $\iiSB{\mathcal{F}}_k(d,f)$,
and the pair $\underbrace{\cdots\iiSB{\mathcal{F}}\iiSC{\mathcal{F}}}_{k
\mbox{ times }}(d,f)$ as $\iiSC{\mathcal{F}}_k(d,f)$, $k = 0,1,2,\ldots$;
with $k=0$ we get the pair $(d,f)$ itself.

(Rigorously speaking, it is necessary to define subsets $S_i$, $i
\in \mathbb{N}$ in the following way: $S_1 = S_\circ \bigcap S_\bullet$,
$S \supset S_1 \supset \ldots \supset S_k \supset$ so as $S_{i+1} =
\{(d,f)\,|\, \iiSC{\mathcal{F}}(d,f) \in S_i\mbox{ and }
\iiSB{\mathcal{F}}(d,f) \in S_i\}$. So, $\iiSC{\mathcal{F}}$ and
$\iiSB{\mathcal{F}}$ are defined on $S_k$ with values in
$S_{k-1}$. Then we have $\iiSC{\mathcal{F}}_k = \underbrace{\ldots
\iiSB{\mathcal{F}}\iiSC{\mathcal{F}}}_{k}$ and $\iiSB{\mathcal{F}}_k =
\underbrace{\ldots \iiSC{\mathcal{F}}\iiSB{\mathcal{F}}}_{k}$ defined
on $S_k$ with values in $S$).

A pair $(d,f)$ is said to be {\it root on $G$}, if for some
$k \in \mathbb{N}_0$ $(d,f) =
\iiSC{\mathcal{F}}_k(g,f_{g})$, if $g \in \iiSB{G}$, or $(d,f) =
\iiSB{\mathcal{F}}_k(g,f_{g})$, if $g \in \iiSC{G}$. Two root pairs
$(d,f_1)$ and $(d,f_2)$ are said to be {\it equivalent}, if $f_1|_{X_d}
\equiv f_2|_{X_d}$. A class of equivalence of the root pair $(d,f)$ we will
denote as $[d,f]$.

\begin{theorem}\label{t_res}

Let $G$~--- connected finite graph without cycles.

1. Following conditions are equivalent:

a) graph $G$ is one of the Dynkin graphs $A_n$, $D_n$, $E_6$, $E_7$,
$E_8$;

b) graph $G$ is finitely representable in $\mathcal{H}$;

c) any finite-dimensional indecomposable representation of a quiver $Q$
is unitarizable.

2. A map $\pi \to (d(i),f(i))$ where $d(i)$ is a dimension of $\pi$ and
$f(i)$ is a character of $\pi$ establishes a one-to-one
correspondence between classes of equivalenc of the indecomposable
locally-scalar representations of a graph $G$ and classes of equivalent
root pairs of a graph $G$.

\end{theorem}


{\bf Proof.}

1. a) implies b). Indeed, if $G$ is a Dynkin graph then all
representations of $G$ are discrete according to the
proposition~\ref{p_dis}. All dimensions of indecomposable representations
are bounded in a whole since there a finite number of them (dimensions of
indecomposable representations are roots of a graph $G$ by the
proposition~\ref{p36}). In order to prove a finite representability of $G$
in $\mathcal{H}$, obviously, it is enough to show that in $\widetilde{\rm
Rep}\,(G,d,f)$ it is contained exactly one, up to equivalence, object
(corollary~\ref{c_equiv}), and that two finite-dimensional representations
$\pi$ and $\widetilde{\pi}$ from ${\rm Rep}\,(G)$ are unitary equivalent if
and only if they have equal dimension $d$, common character $f$ and
$(\pi,f)$ is equivalent to $(\widetilde{\pi},f)$ in ${\rm Rep}\,(G,d,f)$.

b) implies a), since if $G$ is not a Dynkin graph then it is not
finitely representable by the proposition~\ref{p41} (it contains
an extended Dynkin graph).

a) and c) are equivalent by the proposition~\ref{p42}.


2. If two locally-scalar representations $\pi_1$ and $\pi_2$ of a graph
$G$ are unitary equivalent and indecomposable (and, consequently,
finite-dimensional), then it is obvious that their dimensions
are equal, and their characters (uniquely determined) on the common support
are equal too.

Let locally-scalar indecomposable representations $\pi_1$ and $\pi_2$
have equal dimension $d$ and equal character $f$ (on the common
support characters are match, and out of the support we will define
them as equal). Then $\pi_1$ and $\pi_2$ have equal growth $d =
\iiSC{c}_t(\bar{g})$ or $\iiSB{c}_t(\bar{g})$. Thus, pairs $(\pi_1,f)$,
$(\pi_2,f)$ as the objects of the category $\widetilde{\rm
Rep}\,(G,\coprod)$ can be obtained from the simplest object
$(\Pi_{g},f_{g})$ by the functors of the even and odd reflections (see
the proof of the proposition~\ref{p36}) and, therefore, match.


\begin{remark}

Even for Dynkin graphs indecomposable locally-scalar representations are
not determined uniquely by the character: there exist examples
of indecomposable representations with different dimensions but
equal characters.

\end{remark}

\begin{remark}

It is follows from our constructions that for Dynkin graphs (and,
seemingly, only for them) all locally-scalar representations
realize over the field of the real numbers.

\end{remark}


\bigskip

{\Large\it Translated from russian by A.~V.~Roiter.}


\begin{thebibliography}{99}

\bibitem{BerGel} {\it Bernstein~I.N., Gelfand~I.M., Ponomarev~V.A.}
Coxeter functors and Gabriel's theorem. -- IMN, v. XXVIII, pt.~2, p.~19-33
(1973).

\bibitem{Gab} {\it P.~Gabriel.} Unzerlegbare Darstellungen I -- Manuscripta
Math.~6 (1972), 71-103.

\bibitem{Naz} {\it Nazarova~L.A..} Representations of quadriade. -- Izv. AN
USSR, 1967, {\bf 31}, N.~4, p.~1361--1378

\bibitem{GelPon} {\it Gelfand~I.M., Ponomarev~V.A.} Problems of linear
algebra and classification of quadruples of subspaces in a
finite-dimensional vector space. -- Coll. Math. Soc. J. Bolyai, S, Hilbert
Space Operators, Tihany, Hungary, 1970.

\bibitem{Kru} {\it Kruglyak~S.A.} Representations of algebras, which square
of radical eqals zero. Zap. nauch. sem. LOMI, 1972, {\bf 28},
p.~60-69

\bibitem{Naz2} {\it Nazarova~L.A..} Representations of quivers
of infinite type. -- Izv. AN USSR. Ser. mat. 1973, {\bf 37},
p.~752--791.

\bibitem{DonFre} {\it Donovan~P., Freislich~M.R.} The representation theory
of finite graphs and associated algebras. -- Carleton Math. Lecture Notes,
1973, 5, p.1--119.

\bibitem{Kac} {\it Kac~V.G.} Some remarks on representations of quivers and
infinite root systems. -- Lect. Notes Math., 1980, {\bf 832}, p.~311-332.

\bibitem{Kac2} {\it Kac~V.G.} Infinite root systems, representations of
graphs and invariant theory, II. -- I. Algebra, 1982, {\bf 78}, p.~141-162.

\bibitem{GabRoi} {\it P.~Gabriel, A.V.~Roiter} Representations of
Finite-Dimensional Algebras. Springer, 1--171, 1997.

\bibitem{KruSam} {\it Kruglyak~S.A., Samoilenko~Yu.S.} On unitary
equivalence of collections of selfadjoint operators. -- Functional
analysis and its applications, 1980, pt.~1, p.~84-85.

\bibitem{Kru2} {\it Kruglyak~S.A.} Representations of free involutive
quivers. -- In a book Representations and quadratic forms. Kiev: Pub.
of the Inst. of Math. AN UkrSSR, 1979, p.~149-151.

\bibitem{Ser} {\it V.V.~Sergeichuk.} Unitary and Euclidean representations
of a quiver. -- Linear Algebra and its Applications, 278 (1998), p.~37-62.

\bibitem{OstSam} {\it Ostrovskii~V., Samoilenko~Yu.} Introduction to the
Theory of Representations of Finitely Presented $*$-algebras. I.
Representations by bounded operators. Rev. Math. \& Math. Phys., vol.11,
1--261, Gordon and Breach, 1999.

\bibitem{Roi} {\it Roiter~A.V.} Boxes with involution. -- In a book
Representations and quadratic forms. Kiev: Pub.
of the Inst. of Math. AN UkrSSR, 1979, p.~124-126.

\bibitem{RabSam} {\it Rabanovich~V.I., Samoilenko~Yu.S.} When
sum of idempotents or projectors multiple to identity. -- Functional
analysis and its applications, v.~34, pt.~4, p.~91-93 (2000).

\bibitem{KruRab} {\it Kruglyak~S.A., Rabanovich~V.I., Samoilenko~Yu.S.} On
sums of projectors -- Functional analysis and its applications, v.~36,
pt.~3, p.~20-35, 2002.

\bibitem{Kru3} {\it Kruglyak~S.A.} Coxeter functors for one class
of $*$-quivers. -- Ukr. math. jour., v.~54, N.~6, 2002.

\bibitem{Kru4} {\it Kruglyak~S.A.} Coxeter functors for a certain class of
$*$-quivers and $*$-algebras. Methods of Functional Analysis and Topology.
Vol.~8, N.~4, pp. 49--57, 2002.

\bibitem{Bas} {\it H.~Bass.} Algebraic $K$-theory. -- Moscow,
``Mir'', 1973.

\bibitem{RedRoi} {\it I.K.~Redchuk, A.V.~Roiter.} Singular
locally-scalar representation of quivers in Hilbert spaces and
separating functions. -- Present preprint.

\end{thebibliography}
\end{document}